\documentclass[12pt]{article}
\usepackage{fullpage}
\usepackage{multirow}
\usepackage{subcaption}
\usepackage{mathptmx}
\usepackage{relsize}
\usepackage{float}
\usepackage{amsmath}
\usepackage{epsfig}
\usepackage{algorithm}
\usepackage{booktabs}
\usepackage{latexsym}
\usepackage{color}
\usepackage{graphicx}
\usepackage{dsfont}
\usepackage{color}
\usepackage{mathrsfs}
\usepackage{amsmath,amssymb,ifthen,url,graphicx,color,array,theorem}

\usepackage{amssymb,ifthen,url,graphicx,color,array,theorem}

\usepackage{bm}
\usepackage{soul}

\usepackage[flushleft]{threeparttable}
\usepackage{natbib}
 \bibpunct[, ]{(}{)}{,}{a}{}{,}%
\newtheorem{Definition}{\noindent Definition}

\newtheorem{Remark}{Remark}
\newtheorem{Lemma}{Lemma}
\newtheorem{Theorem}{Theorem}

\numberwithin{Theorem}{section}
\numberwithin{Corollary}{section}
\numberwithin{Lemma}{section}
\numberwithin{Proposition}{section}
\numberwithin{Definition}{section}
\numberwithin{Remark}{section}
\numberwithin{algorithm}{section}
\numberwithin{equation}{section}

\newcommand{\mb}[1]{\mathbf{#1}}

\newcommand{\norm}[1]{\left\lVert#1\right\rVert_2}
\newcommand{\abs}[1]{\left|#1\right|}

\makeatletter
\let\save@mathaccent\mathaccent
\newcommand*\if@single[3]{%
  \setbox0\hbox{${\mathaccent"0362{#1}}^H$}%
  \setbox2\hbox{${\mathaccent"0362{\kern0pt#1}}^H$}%
  \ifdim\ht0=\ht2 #3\else #2\fi
  }
\newcommand*\rel@kern[1]{\kern#1\dimexpr\macc@kerna}
\newcommand*\widebar[1]{\@ifnextchar^{{\wide@bar{#1}{0}}}
{\wide@bar{#1}{1}}}
\newcommand*\wide@bar[2]{\if@single{#1}{\wide@bar@{#1}{#2}{1}}
{\wide@bar@{#1}{#2}{2}}}
\newcommand*\wide@bar@[3]{%
  \begingroup
  \def\mathaccent##1##2{%
    \let\mathaccent\save@mathaccent
    \if#32 \let\macc@nucleus\first@char \fi
    \setbox\z@\hbox{$\macc@style{\macc@nucleus}_{}$}%
    \setbox\tw@\hbox{$\macc@style{\macc@nucleus}{}_{}$}%
    \dimen@\wd\tw@
    \advance\dimen@-\wd\z@
    \divide\dimen@ 3
    \@tempdima\wd\tw@
    \advance\@tempdima-\scriptspace
    \divide\@tempdima 10
    \advance\dimen@-\@tempdima
    \ifdim\dimen@>\z@ \dimen@0pt\fi
    \rel@kern{0.6}\kern-\dimen@
    \if#31
      \overline{\rel@kern{-0.6}\kern\dimen@\macc@nucleus\rel@kern{0.4}
      \kern\dimen@}%
      \advance\dimen@0.4\dimexpr\macc@kerna
      \let\final@kern#2%
      \ifdim\dimen@<\z@ \let\final@kern1\fi
      \if\final@kern1 \kern-\dimen@\fi
    \else
      \overline{\rel@kern{-0.6}\kern\dimen@#1}%
    \fi
  }%
  \macc@depth\@ne
  \let\math@bgroup\@empty \let\math@egroup\macc@set@skewchar
  \mathsurround\z@ \frozen@everymath{\mathgroup\macc@group\relax}%
  \macc@set@skewchar\relax
  \let\mathaccentV\macc@nested@a
  \if#31
    \macc@nested@a\relax111{#1}%
  \else
    \def\gobble@till@marker##1\endmarker{}%
    \futurelet\first@char\gobble@till@marker#1\endmarker
    \ifcat\noexpand\first@char A\else
      \def\first@char{}%
    \fi
    \macc@nested@a\relax111{\first@char}%
  \fi
  \endgroup
}
\makeatother

\DeclareMathOperator*{\argmin}{arg\,min}
\DeclareMathOperator*{\argmax}{arg\,max}

\begin{document}

\title{Solving Multi-Objective Optimization via Adaptive Stochastic Search with Domination Measure}
\author{Joshua Q Hale, Helin Zhu, Enlu Zhou\\
H. Milton Stewart School of Industrial and Systems Engineering,\\
Georgia Institute of Technology}
\date{}
\maketitle
\begin{abstract}
For general multi-objective optimization problems, we propose a novel performance metric called domination measure to measure the quality of a solution, which can be intuitively interpreted as the size of the portion of the solution space that dominates that solution. As a result, we reformulate the original multi-objective problem into a stochastic single-objective one and propose a model-based approach to solve it.
We show that an ideal version algorithm of the proposed approach converges to a set representation of the global optima of the reformulated problem.
We also investigate the numerical performance of an implementable version algorithm by comparing it with numerous existing multi-objective optimization methods on popular benchmark test functions.
The numerical results show that the proposed approach is effective in generating a finite and uniformly spread approximation of the Pareto optimal set of the original multi-objective problem, and is competitive to the tested existing methods.
\end{abstract}

\section{Introduction}\label{Sec1:Introduction}

Problems that require optimizing several objectives concurrently are known as multi-objective optimization problems.
Obviously, this type of problems arises in many real-world applications, including construction science (\cite{le2003coupled}), economics (\cite{toffolo2002evolutionary}), medical treatments (\cite{qasem2011radial}), and logistics (\cite{lee2008multi}), in which incommensurable and conflicting objectives need to be optimized. Therefore, it is often unlikely to have a solution that optimizes all objectives simultaneously.
A more reasonable goal is to obtain a set of solutions, where the quality of each solution is incomparable without any prior knowledge of preference.
These solutions are known as Pareto optimal solutions, which are ``optimal'' in the sense that no other solutions in the solution space are superior to them while taking into account all of the objectives.
The set of Pareto optimal solutions is called the Pareto optimal set, and its image in the objective space is called the Pareto front.

Numerous methods have been developed to find or approximate the Pareto optimal set of a general multi-objective optimization problem,
among which perhaps the most popular ones are evolutionary algorithms (\cite{deb2001multi}) that use iterative selection, mutation and crossover operations to generate multiple Pareto optimal solutions in parallel.
Other methods include stochastic search methods (\cite{zabinsky2013}) that choose candidate solutions at random and improve the way those candidate solutions are selected in each iteration, and particle swarm methods (\cite{mete2014multiobjective}) that keep a population of potential solutions (particles) that are manipulated by a velocity vector, which changes the position of the particles at each iteration.
Furthermore, simulated annealing (\cite{czyzzak1998}) has also been implemented in the multi-objective domain by using multiple weight vectors to convert the problem into several single-objective problems.

The aforementioned methods can be classified as instance-based methods, as they maintain a population of candidate solutions, and new candidate solutions are generated by only considering previously generated solutions.
In contrast, model-based search algorithms generate candidate solutions from parameterized sampling distributions that are iteratively updated using previous candidate solutions.
They are effective at optimizing functions that lack structural properties such as differentiability and convexity.
Although model-based algorithms have mostly been used to solve single-objective optimization problems, there are some methods that incorporate the Cross Entropy (CE) (\cite{rubinstein2001combinatorial}) method, a model-based approach, to solve multi-objective problems, see e.g., \cite{unveren2007multi} and \cite{bekker2011cross}.
In \cite{unveren2007multi} samples are generated from a fixed number of sampling distributions, which allows the algorithm to be very fast and simple to implement. However, since a fixed number of sampling distributions are used, this method may have difficulty with constructing a sampling distribution that has majority of its mass on isolated points of the Pareto optimal set. \cite{bekker2011cross} propose a histogram approach in which candidate solutions are drawn independently along each dimension of the solution space.
Although this method performs very well on some problems, it could have difficulty in capturing the entire Pareto optimal set that consists of highly correlated solutions.

Most of the methods mentioned above are designed for optimizing multiple objectives that are deterministic. There also exist a number of methods to solve the stochastic version of this problem, see, e.g, \cite{li2015mo} and \cite{2015felhunpasWSC}. The details of these methods are beyond the scope of this paper since we are only focused on deterministic objectives.

In view of the existing methods in the literature, a method that explores the superiority or the dominance relationship among solutions by simple quality metrics is still lacking.
Incorporating such metrics reduces problem dimensionality since the multi-dimensional objective space is mapped onto a single-dimensional one through direct comparisons among solutions.
Thus, the original multi-objective problem is transformed into a single-objective one, in which the objective is to find the solutions that optimize a particular quality metric.
Of course, as pointed out by \cite{zitzler2003performance}, a reduction of problem dimensionality will inevitably cause the loss of information.
That is, generally the global optimal set of the reformulated problem will not be the same as the Pareto optimal set.
In particular, \cite{zitzler2003performance} prove that in order for a metric to retain the Pareto dominance relation the dimension of the metric and the objective space must be equivalent.
There are some initial explorations along this line of research.
For example, a quality metric can be derived by a weighted aggregation of the objective functions.
The issue with this quality metric is that a single choice of weights leads to at most one point in the Pareto optimal set.
Although multiple points could be obtained via changing the weights, the final approximation could be unsatisfactory since a uniform spread of the weighting coefficients does not necessarily produce a uniform spread of Pareto optimal solutions.
Moreover, this technique is infeasible for problems with a large amount of objectives, since the total number of weight combinations grows exponentially with respect to (w.r.t) the number of objectives.

\cite{zitzler2003performance} proposed a quality metric of a solution set termed hypervolume, which is defined as the total size of the area that is dominated by that set in the objective space w.r.t. a reference point.
It is desirable for performance assessment of a solution set since it has the ability to measure how close solutions are to the Pareto front as well as how evenly spread the solutions are in the Pareto optimal set.
Therefore, it can qualitatively compare two different solution sets by the use of one value.
Since the hypervolume metric is influenced by any type of progress to the true Pareto front, many methods have used this metric as a way to guide the search to favorable areas of the solution space, see \cite{fleischer2003measure}, etc.
Although the hypervolume metric has many favorable properties,  there also are some drawbacks.
In particular, the choice of the reference point could affect progression of the search to promising areas of the solution space, the hypervolume is biased towards convex regions of the objective space, and as the number of objectives increase so does the computational complexity of calculating the hypervolume.

In this paper, we introduce a new parameter-free and unary performance metric to measure the quality of solutions, termed as \textit{domination measure}.
Practically, the domination measure of a solution can be viewed as the measure of the region in the solution space that dominates that solution.
Unlike hypervolume mentioned previously, it is a global performance metric since all the solutions in the solution space are taken into account for its computation.
To the best of our knowledge, this work is among the very first to incorporate such a performance metric for solving multi-objective problems.

There are many advantages for using domination measure as a performance metric.
Foremost, unlike the aforementioned scalarization methods, this metric does not require any tuning of parameters.
Moreover, the domination measure of a solution is a rigorous quantification of the quality of a solution. That is, the lower the domination measure, the better the solution. Finally, if a solution is Pareto optimal, then it has a domination measure of zero since no solution dominates it. When the solution space consists of finite solutions, the set of Pareto optimal solutions is exactly the set of solutions with domination measure of zero.
Precisely, a solution being Pareto optimal is equivalent to the solution having a domination measure of zero.
In contrast, for a continuous solution space a solution that is not Pareto optimal can also have a domination measure of zero, which is not surprising since some of the dominance relation is lost by reducing the dimensionality of the objective space through domination measure. 
Although theoretically the set of solutions with domination of zero is not equivalent to the Pareto optimal set, it is a sufficiently good approximation in the sense that no solution dominates it
almost surely.

By employing domination measure as a performance metric, we are able to transform the original multi-objective optimization problem into a single-objective problem, where the goal is to find solutions that have a domination measure of zero. Note that this objective is also stochastic since domination measure is an expectation of an indicator function on the dominance relation w.r.t. a uniform probability measure.
While we benefit from a significant reduction in problem complexity, we make the compromise from finding Pareto optimal solutions to finding solutions with a domination measure of zero. Nevertheless, we will propose an approach tailored for the reformulated problem and show empirically that it is not prone to converge to a non-Pareto optimal solution with a domination measure of zero.

To solve the reformulated stochastic single-objective problem, we propose a model-based approach that finds multiple global optimal solutions. The idea is to introduce a mixed sampling distributions with components from a parameterized family of densities, and iteratively update the sampling distribution parameters by minimizing the Kullback-Leibler divergence between properly constructed reference distributions and the parameterized sampling distributions. The sampling distributions are updated until they are degenerate  distributions concentrated on a global optima. Our approach is similar to existing model-based approaches such as the CE method and Model Adaptive Reference Search (MRAS) (\cite{hu2007model}) in the sense that reference distributions are introduced to guide the sampling process while the actual sampling is achieved by the parameterized distributions. A major difference from the aforementioned methods is that in every iteration our approach keeps track of a population of sampling distributions with an adaptive number of components instead of a single sampling distribution. This ensures that multiple uniformly spread global optimal solutions are generated rather than a single optimal solution.

Based on the proposed approach we design two algorithms, which are an ideal version and a implementable version. For the ideal version algorithm, we show that for every solution that has a domination measure of zero, there exists a sequence of parameterized sampling distributions that converges to the degenerate distribution on that solution.
For the implementable version algorithm, we show empirically that it is competitive to many existing multi-objective optimization methods by testing them on several classic benchmark problems.
In particular, we observe that in all the cases tested, our approach is able to generate solutions that are approximately uniformly spread across the Pareto optimal set by a user-specified threshold distance.

In summary, the contributions in this paper are as follows:
\begin{itemize}

\item For multi-objective optimization problems, we introduce a novel performance metric termed domination measure to determine the quality of a solution, and we reformulate the original problem into a stochastic single-objective one that seeks to minimize the domination measure. 

\item We propose a novel model-based approach to solve the reformulated problem, leading to an ideal version algorithm that possesses nice convergence properties and an implementable version algorithm that performs well numerically.

\item We show that our proposed approach produces a finite and uniformly spread approximation of the Pareto optimal set and performs competitively to or even outperforms many existing approaches. 
\end{itemize}

The remainder of this paper is organized as follows:
In Section \ref{Sec2:Formulation}, we introduce the concept of domination measure under the setting of multi-objective optimization, and use it to reformulate the original problem into a stochastic single-objective one.
In Section \ref{Sec3:Model}, we present a model-based approach to solve the reformulated problem, and describe the ideal version and implementable version algorithms.
In Section \ref{Sec4:Numerical}, we conduct numerical experiments to demonstrate the effectiveness and advantages of the proposed approach by comparing it with existing approaches. Finally, conclusions and future directions of research are given in Section \ref{Sec5:Conclusions}.

\section{Multi-Objective Optimization and Domination Measure}\label{Sec2:Formulation}
A general multi-objective problem consists of minimizing (or maximizing) multiple objectives over a defined solution space, which can be formulated as follows:
\begin{equation}\label{eq.2.1}
\begin{aligned}
&\min \; &&\mb{f}(x)=\{f_1(x), f_2(x),...,f_n(x)\}\\
&\text{s.t.}\; &&x\in \mathcal{X},
\end{aligned}
\end{equation}
where $\mathcal{X}$ denotes the solution space that might be described by constraints, and $\{f_i(\cdot): \mathcal{X} \rightarrow \mathbb{R},\; i=1,...,n\}$ are scalar functions. Without loss of generality, we assume $\mathcal{X}$ is a bounded subset of $\mathbb{R}^d$ and all objective are minimized.

Since it is rarely the case that there exists a solution that minimizes all the objectives simultaneously, a reasonably compromised goal is to find all the solutions that are not dominated by any other solutions in $\mathcal{X}$ in terms of the objective values.
Specifically, a solution $x\in\mathcal{X}$ is (Pareto) dominated by another solution $y\in\mathcal{X}$ if $f_i(y)\le f_i(x)$ for all $i=1,...,n$, and there exists one $j\in\{1,...,n\}$ such that $f_j(y)<f_j(x)$.
In other words, a solution $x$ is dominated by another solution $y$ if all the objectives evaluated at $y$ are better than (less than or equal to) the ones evaluated at $x$, and at least one objective evaluated at $y$ is strictly better than (less than) the one evaluated at $x$.
Note that Pareto dominance is a (strict) partial order defined on $\mathcal{X}$ since it is irreflexive, asymmetric and transitive. For simplicity, we use $y \prec_{d} x$ to denote that $y$ dominates $x$.
Following this, a solution $x\in\mathcal{X}$ is called Pareto optimal if all other solutions in $\mathcal{X}$ do not dominate $x$, and thus the goal is to find the Pareto optimal set or a uniformly spread subset representation of the Pareto optimal set.

In general, this problem is difficult because the dominance relation is a partial order defined on the solution space, and the problem is essentially a combinatorial problem over a solution space that is often continuous.
As mentioned before, many approaches in the literature are developed for solving various approximations or reductions of a multi-objective optimization problem.
In particular, one type of approach is to reformulate the original problem into a single-objective one (e.g., through a weighting scheme on the objective functions) and apply algorithms that are designed for problems with a single objective.
As a result of the reformulation, some information of the original problem (e.g., the dominance relationship or the defined partial order among solutions) is lost.

In the next subsection, we will propose a performance metric called domination measure that quantifies the dominance relationship between a solution of interest and all other solutions in the solution space, and use it to reformulate the original multi-objective problem into a stochastic single-objective one.
In contrast to most of the existing reformulation techniques, the proposed reformulation technique relies on the combinatorial dominance relationship between solutions rather than the absolute objective values of the solutions.
We will show that it leads to a close approximation of the original problem, and is solvable by a proposed model-based optimization method.

\subsection{Domination Measure}\label{Sec2.subsec1:Domination}

Roughly speaking, domination measure of a solution describes the portion of the solutions in the solution space that dominates that solution.
To ease presentation, let us consider defining domination measure of a solution $x\in \mathcal{X}$, denoted by $D(x)$, in the following two cases: 1) the solution space $\mathcal{X}\in \mathbb{R}^d$ has a non-zero and finite measure w.r.t. the Lebesgue measure $\nu$ of the $d$-dimensional Euclidean space $\mathbb{R}^d$; 2) $\mathcal{X}$ is finite.
For the second case, we extend the definition of Lebesgue measure on $\mathbb{R}^d$ to any finite set $S$ with cardinality $|S|$ as follows:
Suppose $S_1\subseteq S$ is a subset of $S$, then we define the Lebesgue measure of $S_1$ by $\nu(S_1)=|S_1|/|S|$, i.e., $\nu(S_1)$ is the cardinality of $S_1$ divided by the cardinality of $|S|$; in particular, $\nu(S)=1$.
We formally define domination measure as follows:

\begin{Definition}\label{def.2.1}
\textbf{(Domination Measure)}. Assume $\mathcal{X}$ in problem (\ref{eq.2.1}) is either a subset of $\mathbb{R}^d$ that has a non-zero and finite Lebesgue measure or a finite set.
Further assume for all $x\in \mathcal{X}$ the set of solutions that dominates $x$, denoted by $\mathcal{D}_x$, is Lebesgue measurable.
Then, for all $x\in \mathcal{X}$, the domination measure $D(x)$ of $x$ is defined as
\begin{equation}\label{eq.2.2}
D(x)\overset{\triangle}=\frac{\nu(\mathcal{D}_x)}{\nu(\mathcal{X})},
\end{equation}
where $\nu(A)$ is the Lebesgue measure of $A$.
\end{Definition}

\begin{Remark}\label{rem.2.1}
It is easy to extend Definition \ref{def.2.1} to the case where the set $\mathcal{X}$ is (countable or uncountable) infinite and has Lebesgue measure of zero. Specifically, let $D(x)=1$ if the cardinality of $D_x$ is infinite and $0$ otherwise. Since for this case the domination measure of a solution only has two trivial values ($0$ and $1$), we omit the discussion about this case.
\end{Remark}

Intuitively, the domination measure $D(x)$ of $x$ is the ratio of the measure of $\mathcal{D}_x$  to the measure of the entire solution space $\mathcal{X}$.
Therefore, $\forall x\in \mathcal{X}$ we have $0\le D(x)\le 1$. Furthermore, it is a rigorous performance metric for the quality of a solution. If the domination measure of a solution is close to zero, then that solution is dominated by a small number of solutions in the solution space.

\begin{Lemma}\label{lem.2.1}
For any Pareto optimal solution $x^\ast \in \mathcal{X}$, its domination measure $D(x^\ast)=0$.
\end{Lemma}

\begin{figure}[htb!]
\centering
\includegraphics[width=0.5\linewidth]{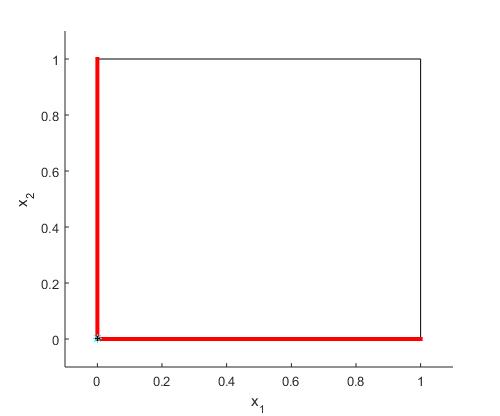}
\caption{Illustration of Domination Measure.}
\label{fig.2.1}
\end{figure}

The other direction of the statement in Lemma \ref{lem.2.1} might not be true, which means a solution with domination measure of zero might not be Pareto optimal.
To illustrate this, consider a simple example with two objectives and a two-dimensional solution space $[0,1]\times[0,1] \in \mathbb{R}^2$.
The two objectives are $f_1(x_1, x_2)=x_1$ and $f_2(x_1, x_2)=x_2$. Therefore, $(f_1, f_2)$ is a mapping from $[0,1]\times[0,1]$ to itself.
The solution space for this example is displayed in Figure \ref{fig.2.1}.
Although $(0,0)$ is the lone Pareto optimal solution, all solutions in red have a domination measure of zero.
For example, point $(0.4,0)$ is dominated by all points where $x \in [0,0.4)$ and $y=0$, but the measure of those points is equal to zero. Conversely, for a finite solution space, no information is lost using domination measure. 

Since a Pareto optimal solution achieves a minimum domination measure of zero, domination measure retains part of the information from the partial order dominance relation. The benefit is that the original multi-dimensional combinatorial problem can be reduced to a simple single-objective one with a known minimum objective value of zero.

\begin{figure}
\centering
\begin{subfigure}{.5\textwidth}
  \centering
  \includegraphics[width=1\linewidth]{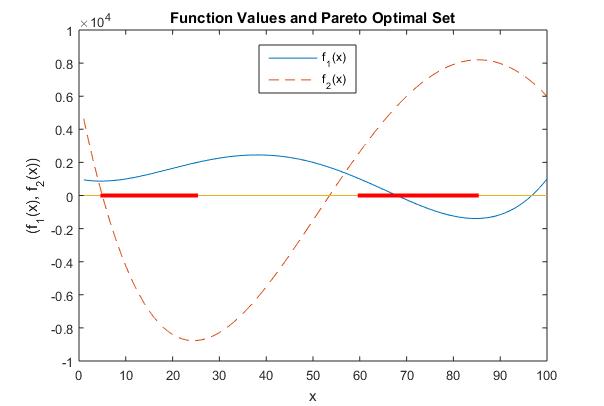}
  \caption{Pareto Optimal Set by Function Values.}
  \label{fig.2.2.sub1}
\end{subfigure}%
\begin{subfigure}{.5\textwidth}
  \centering
  \includegraphics[width=1\linewidth]{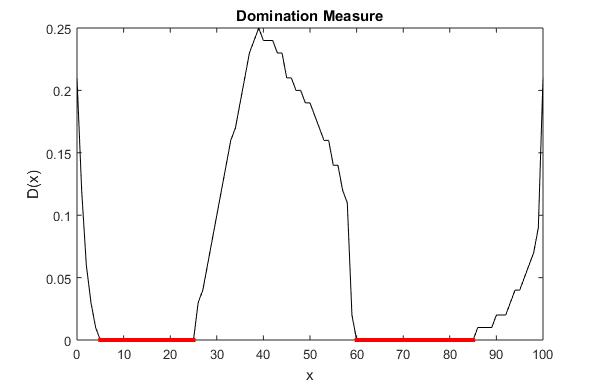}
  \caption{Pareto Optimal Set by Domination Measure}
  \label{fig.2.2.sub2}
\end{subfigure}
\caption{Pareto Optimal Set by Function Values and Domination Measure.}\label{fig.2.2}
\end{figure}

To demonstrate the effectiveness of this reduction, we consider an example from \cite{mete2014multiobjective} as follows:
The solution space is $\mathbb{Z}\cap[0,100]$, i.e., the set of all the integers between $0$ and $100$, which is a finite set.
There are two objective functions:
$$
f_1(x) = 0.001x(x-10)(x-60)(x-100) + 1000,\;
f_2(x) = 0.001x(x-70)(x-100)(x-200) + 6000.
$$
The Pareto optimal set is $P^* = \mathbb{Z} \cap ([5,25] \cup [60,85])$, and it is highlighted in Figure \ref{fig.2.2}.
Locating the Pareto optimal set directly using function values (see Figure \ref{fig.2.2.sub1}) is difficult, whereas the Pareto optimal set is easily identifiable using domination measure (see Figure \ref{fig.2.2.sub2}). That is, the solutions with domination measure zero are the Pareto optimal set for a finite solution space.

We minimize the domination measure instead of solving the original combinatorial problem.
Note that we can reformulate (\ref{eq.2.2}) as
\begin{equation}\label{eq.2.3}
D(x)=\frac{\nu(\mathcal{D}_x)}{\nu(\mathcal{X})}
=\frac{\int_{\mathcal{X}}\mathds{1}\left\{y\prec_d x\right\}\nu(dy)}{\int_{\mathcal{X}}\nu(dy)}=
\int_{\mathcal{X}}\mathds{1}\left\{y\prec_d x\right\}U_\mathcal{X}\nu(dy)=
\mathbb{E}_U\left[\mathds{1}\left\{y\prec_d x\right\}\right],
\end{equation}
where $\mathds{1}\left\{E\right\}=1$ if the event $E$ is true and $\mathds{1}\left\{E\right\}=0$ otherwise, $U_\mathcal{X}(\cdot)$ is the uniform probability measure (induced by the uniform distribution) on $\mathcal{X}$, and $\mathbb{E}_U[\cdot]$ denotes the expectation w.r.t. $U_\mathcal{X}(\cdot)$.
Consequently, we solve the following single-objective stochastic optimization problem:
\begin{equation}\label{eq.2.4}
\min_{x\in\mathcal{X}} \;D(x)=\mathbb{E}_U\left[\mathds{1}\left\{y\prec_d x\right\}\right], \quad \mbox{or equivalently} \quad \max_{x\in\mathcal{X}} \;-D(x)=\mathbb{E}_U\left[-\mathds{1}\left\{y\prec_d x\right\}\right].
\end{equation}
Typically, the goal is to find multiple global optima of problem (\ref{eq.2.4}) such that they approximately form a uniform distribution on the Pareto optimal set of the original problem (\ref{eq.2.1}).
Suppose the Pareto optimal set is nonempty, then by Lemma \ref{lem.2.1}, there always exists an $x\in\mathcal{X}$ such that $D(x)=0$. As a result, problem (\ref{eq.2.4}) always admits at least one global optimal. For simplicity, we denote the set of global optima of problem (\ref{eq.2.4}) by $A^\ast$. Therefore, $\forall x\in\mathcal{X}$, $D(x)=0$ if and only if $x\in A^\ast$.

Explicitly calculating $D(x)$ is usually computationally expensive or infeasible, especially in the case of a continuous solution space.
Alternatively, since Monte Carlo simulation is straightforward to implement and scales well with the dimension of the solution space, it is a good choice for estimating $D(x)$. The general procedure is as follows:
Draw independent and identically distributed (i.i.d.) samples according to the uniform distribution $U_\mathcal{X}(\cdot)$ on the solution space, which is denoted by $\{x^1,..., x^N\}$. Then calucalte \begin{equation}\label{eq.2.5}
\widetilde{D}(x^i)\overset{\triangle}=
\frac{1}{N}\sum_{j=1}^N \mathds{1}\left\{x^j\prec_d x^i\right\},
\end{equation}
which is an unbiased estimator of $D(x^i)$.
We also point out that the sampling distribution could be of other forms via the principle of importance sampling as long as it is fully supported on $\mathcal{X}$.

\section{A Framework of Model-based Approach}\label{Sec3:Model}

Recall that we aim to find multiple global minimizers of the domination measure $D(x)$ such that they approximately form a uniform distribution on the Pareto optimal set of problem (\ref{eq.2.1}).
Since $D(x)$ is often evaluated via simulation, its structural properties such as convexity or differentiability are unknown.
Thus, traditional gradient-based optimization methods might not be applicable for the minimization of $D(x)$.
In contrast, model-based optimization methods are good alternatives as they impose minimal requirements on the problem structure.
Common model-based methods include Annealing Adaptive Search (ASS) (\cite{romeijn1994simulated}), the Cross-Entropy (CE) method (\cite{rubinstein2001combinatorial}), Model Reference Adaptive Search (MRAS) (\cite{hu2007model}), and Gradient-based Adaptive Stochastic Search (GASS) (\cite{zhou2014gradient}), etc.

The main idea of model-based methods is to introduce a sampling distribution, which often belongs to a parameterized family of densities, over the solution space, and iteratively update the parameters of the sampling distribution by generating and evaluating candidate solutions. Specifically, the methods iteratively carry out the following two steps:

\begin{itemize}
\item[1.] Generate candidate solutions according to the sampling distribution.
\item[2.] Based on the evaluation of the candidate solutions, update the parameter of the sampling distribution.
\end{itemize}

The hope is to have the sampling distribution more and more concentrated on the promising region of the solution space where the optimal solutions are located, and eventually become a degenerate distribution on one of the global optima.
Therefore, finding an optimal solution in the solution space is transformed to finding an optimal sampling distribution parameter in the parameter space.
A key difference among the aforementioned model-based methods lies in how to update the sampling distribution parameter.
For example, in CE and MRAS the updating rule is derived by minimizing the Kullback-Leibler (KL) divergences between a converging sequence of reference distributions and a chosen exponential family of densities.

Compared with gradient-based methods, model-based methods are more robust in the sense that at every iteration they exploit the promising region of the solution space that has already been identified, while maintaining exploration of the entire solution space.
The updating rule on the sampling distribution parameter controls the balance between exploration and exploitation.

In principle, we could extend all the aforementioned model-based methods to our problem setting in a direct manner. However, there are two main issues with implementing existing model-based methods:
1) all these methods are designed for producing a single global optimal solution which corresponds to a degenerate parameterized distribution, but multiple evenly spread global optimal solutions are needed to represent the Pareto optimal set;
2) the methods that are convergent require the uniqueness of the global optimal solution, which in general is not satisfied for problem (\ref{eq.2.4}).

To address the issues mentioned above, we propose a model-based method with a mixed sampling distribution that consists of a number of distributions from the same parameterized family of densities, and iteratively update the parameters of the mixed sampling distribution. The hope is that each component of the mixed sampling distribution will concentrate on a promising region of the solution space and explore that region exclusively.
Eventually, all the components of the mixed sampling distribution will be become degenerate distributions concentrated on distinct global optimal solutions.

The choice of the number of components in the mixed sampling distribution and  the method for updating each component are important.
Since the promising region of the solution space is unknown and can only be explored via sample evaluations, the number of components needs to be determined adaptively so that the all the areas containing the global optima are eventually explored.
The updating scheme on the components needs to achieve convergence so that degenerate sampling distributions on individual global optimal solutions are obtained.

In essence, our updating rule on the sampling distributions is similar to the one in the CE or MRAS, in which the sampling distribution parameter at each iteration is derived by minimizing the KL divergence between a specific reference distribution and a parameterized family of densities.
However, our updating rule keeps track of an increasing population of reference distributions instead of a single one at each iteration. A common choice of the parameterized family of densities is the exponential family defined as follows:

\begin{Definition}\label{def.3.1} \textbf{Exponential Famliy}.
A parameterized family $\{g(x;\theta): \theta\in \Theta\}$ is an exponential family of densities if it satisfies
\begin{equation}\label{eq.3.1}
g(x;\theta)=\exp\left\{\theta^T \Gamma(x)-\eta(\theta)\right\},
\end{equation}
where $\Gamma(x)=[\Gamma_1(x),...,\Gamma_{d_\theta}(x)]^T$ is the vector of sufficient statistics, $\eta(\theta)=\ln\{\int\exp(\theta^T \Gamma(x))dx\}$ is the normalization factor to ensure $g(x;\theta)$ is a probability density function, and $\Theta=\{\theta: |\eta(\theta)|<\infty\}$ is the natural parameter space with a nonempty interior. We assume that $\Gamma(\cdot)$ is a continuous mapping.
\end{Definition}

We point out that many common probability distributions belong to the exponential family, including Gaussian, Poisson, Binomial, Geometric, etc. Within the exponential family of densities, the updating of the sampling parameter using a reference sampling distribution is carried out as follows:
Suppose $h(x)$ is a reference sampling distribution of interest, then the corresponding sampling distribution $g(x;\theta_\ast)$ within the exponential family of densities could be found via

where $\mathbb{E}_h[\cdot]$ denotes the expectation w.r.t. $h(\cdot)$. It is important to note that the problem (\ref{eq.3.2}), i.e., the minimization of KL divergence, can be carried out in analytical form since $g(\cdot;\theta)$ belongs to the exponential family.
Specifically, by Definition \ref{def.3.1}, problem (\ref{eq.3.2}) is equivalent to
\begin{equation*}
\argmax_{\theta\in \Theta} \int_{x\in \mathcal{X}}
(\theta^T \Gamma(x)-\eta(\theta))h(x)dx.
\end{equation*}
Note that $(\theta^T \Gamma(x)-\eta(\theta))$ is strictly concave in $\theta$ (see, e.g., \cite{lehmann2006theory}).
It follows that
$\int_{x\in \mathcal{X}}
(\theta^T \Gamma(x)-\eta(\theta))h(x)dx$
is also strictly concave in $\theta$.
Therefore, (\ref{eq.3.2}) admits a unique optimal solution $\theta_\ast$ that satisfies the first-order condition as follows:
\begin{equation*}
\int \left(\Gamma_j(x)-\frac{\int \Gamma_j(x) \exp
\left(\theta_\ast^T\Gamma(x)\right)dx}{\int \exp
\left(\theta_\ast^T\Gamma(x)\right)dx}\right)h(x)dx=0,~~~ j=1,...,d_\theta,
\end{equation*}
or equivalently,
\begin{equation}\label{eq.3.3}
\mathbb{E}_h[\Gamma(x)]=\mathbb{E}_{\theta_\ast}[\Gamma(x)],
\end{equation}
where $ \forall \theta \in \Theta,~ \mathbb{E}_\theta[\cdot]$ denotes the expectation w.r.t. $g(x;\theta)$.

In practice, finding the exact value of $\theta_\ast$ can be difficult since $\mathbb{E}_h[\Gamma(x)]$ usually does not admit a closed-form expression;
however, an good estimate of $\theta_\ast$ could be obtained via the principle of importance sampling, noting that
\begin{equation*}
\mathbb{E}_{\theta_\ast}[\Gamma(x)]=\mathbb{E}_h[\Gamma(x)]=
\mathbb{E}_{\theta^\prime}
\left[\frac{\Gamma(x)h(x)}{g(x;\theta^\prime)}\right],
\end{equation*}
and i.i.d. samples can be drawn according to $g(\cdot;\theta^\prime)$ to estimate the right hand side of above equation.

For example, if $\{g(\cdot;x)\}$ is the family of multivariate Gaussian distributions $\mathcal{N}(\mu,\Sigma)$, where $\mu$ is the mean parameter, $\Sigma$ is the covariance matrix, and $\theta=(\mu, \Sigma)$. Then $\theta_\ast=(\mu_\ast, \Sigma_\ast)$ can be solved as
\begin{equation*}
\mu_\ast=\frac{\mathbb{E}_{\theta^\prime}
\left[h(x)x/g(x;\theta^\prime)\right]}{\mathbb{E}_{\theta^\prime}
\left[h(x)/g(x;\theta^\prime)\right]}
~~\text{and}~~
\Sigma_\ast=\frac{\mathbb{E}_{\theta^\prime}
\left[h(x)(x-\mu_\ast)(x-\mu_\ast)^T/g(x;\theta^\prime)\right]}
{\mathbb{E}_{\theta^\prime}
\left[h(x)/g(x;\theta^\prime)\right]},
\end{equation*}
and estimated via sampling according to $g(\cdot;\theta^\prime)$.

\subsection{Algorithm---An Ideal Version}
\label{sec3.subsec1:Ideal}

We first present an ideal version algorithm of the proposed approach.
Although it is not implementable in practice, its merit lies in revealing the mathematical intuition of an implementable version algorithm, which will be introduced later.

Denote the population of reference sampling distribution at the $k^{th}$ iteration by $\{h_{i,k}(x): i=1,...,I_k\}$, where $I_k$ is the number of sampling distributions in the population.
Let $\{\gamma_k: k=1,...\}$ be a decreasing sequence of reference values that satisfies $\gamma_k \in (0,1]$ and $\lim_{k\rightarrow \infty}\gamma_k=0$.
We will use $\{\gamma_k\}$ as reference values for characterizing elite solutions and promising regions of the solution space in terms of the domination measure.
Specifically, let $A_k :=\{x\in\mathcal{X}: D(x)\le \gamma_k\}$ be the set of solutions with domination measure value below $\gamma_k$. Loosely speaking, $A_k$ will approach $A^\ast$ as $k\rightarrow \infty$, $A^\ast$ is the set of solutions with domination measure $0$.

Let us further introduce a partition on the set $A_k$, denoted by $\pi_k:=\{A_{1,k},...,A_{I_k,k}\}$, such that each $A_{i,k}$ is $U$-measurable (recall that $U$ is the uniform measure on $\mathcal{X}$), and
\begin{equation*}
A_k=\cup_{i=1}^{I_k} A_{i,k}
~~\text{and} ~~
A_{i,k} \cap A_{j,k}=\emptyset, \forall i\neq j.
\end{equation*}
Note that the cardinality of $\pi_k$ is equal to $I_k$, which is the number of the sampling distributions at the $k^{th}$ iteration. Define the magnitude $\lvert \pi_k\rvert$ of the partition $\pi_k$ by
$
\lvert \pi_k\rvert := \max\limits_i diam(A_{i,k}),
$
where $diam(A):=\sup_{x,y}\norm{x-y}$ is the ``diameter'' of set $A$ and $\norm{\cdot}$ is the vector Euclidean norm. We impose a requirement on $\{\pi_k: k=1,...\}$ as follows:
\begin{equation}\label{eq.3.4}
\lim_{k\rightarrow \infty} \lvert \pi_k \rvert =0.
\end{equation}
Requirement \ref{eq.3.4} forces the partition to become arbitrarily fine for sufficiently large $k$'s.
Equivalently, each partitioning subset $A_{i,k}$ is shrinking to a degenerate point in $A^\ast$ or a empty set as $k\rightarrow \infty$. Such a partition always exists, since $\mathcal{X}$ is bounded

Let the reference distribution $h_{i,k}(x)\propto \mathds{1}\{x\in A_{i,k}\}$, where $h_{i,k}(x)$ is the uniform distribution supported on $A_{i,k}$. Denote the corresponding sampling distribution by $g(x;\theta_{i,k})$, where
\begin{equation*}
\theta_{i,k}=\argmin_{\theta\in\Theta} KL(h_{i,k}(\cdot),
g(\cdot;\theta)).
\end{equation*}
Note that the way to construct the reference sampling distribution is similar to the one in the CE method (see \cite{rubinstein2001combinatorial}). Of course, one could also use more sophisticated reference sampling distributions by an introduction of a shape function to put non-equal weights for $x\in A_{i,k}$ (see \cite{hu2007model}).
The mixed sampling distribution at the $k^{th}$ iteration, denoted by $g_{k}(x)$, consists of $g(x;\theta_{i,k})$ with equal weights. That is,
\begin{equation}\label{eq.3.5}
g_k(x)=\frac{1}{I_k}\sum_{i=1}^{I_k} g(x;\theta_{i,k}).
\end{equation}
The complete algorithm is summarized in the following Algorithm \ref{alg.3.1}, which is also referred to as Algorithm ``$\mathbf{SASMO_0}$''.

\begin{algorithm}[htb]
\caption{\textbf{S}tochastic \textbf{A}daptive \textbf{S}earch for \textbf{M}ulti-objective \textbf{O}ptimization---Ideal}
\label{alg.3.1}
1. \textbf{Initialization}:  Choose a parameterized family of densities $\{g(x;\theta): \theta\in \Theta\}$. Specify a sequence of reference values $\{\gamma_k: k=1,...\}$.
\\
2. \textbf{Iteration}: For the $k^{th}$ iteration, choose a partition $\pi_k=\{A_{1,k},...,A_{I_k,k}\}$ on $A_k$, where $A_k=\{x\in\mathcal{X}: D(x)\le \gamma_k\}$. Then determine the sampling distribution parameters $\{\theta_{i,k}\}$ by
\begin{equation*}
\theta_{i,k}=\argmin_{\theta\in\Theta} KL(h_{i,k}(\cdot),
g(\cdot;\theta)).
\end{equation*}
3. \textbf{Termination}: Check if some stopping criterion is satisfied. If yes, stop and return the means of the currents sampling distributions; else, set $k:=k+1$ and go back to step 2.
\end{algorithm}

Let us analyze the convergence properties of Algorithm ``$\mathbf{SASMO_0}$'', which provides an intuitive understanding towards the convergence of its implementable version.
In particular, we will show that the parameterized sampling distributions $\{g(x;\theta_{i,k}\}$ converge to degenerate distributions on the global optima of the reformulated problem (\ref{eq.2.4}) if $g(x;\theta)$ belongs to an exponential family of densities, as summarized in the following Theorem \ref{thm.3.1}.

\begin{Theorem}\label{thm.3.1}
Suppose the parameterized family of densities $\{g(x;\theta): \theta\in \Theta\}$ is an exponential family defined by Definition \ref{def.3.1}.
Further suppose that the sequence of reference values $\{\gamma_k\in (0,1]: k=1,...\}$ and the sequence of partitions $\{\pi_k: k=1,...\}$ satisfy, respectively,
\begin{equation}\label{eq.3.6}
\lim_{k\rightarrow \infty} \gamma_k=0 ~~\text{and}~~
\lim_{k\rightarrow \infty} \lvert \pi_k \rvert =0.
\end{equation}
Then $\forall x^\ast \in A^\ast$, there exists a sequence of sampling distributions $\{g(x;\theta_{i_k, k}): k=1,...\}$, where $i_k \in \{1,...,I_k\}$, such that
\begin{equation}\label{eq.3.7}
\lim_{k\rightarrow\infty}\mathbb{E}_{\theta_{i_k, k}}
\left[\Gamma(x)\right]=\Gamma(x^\ast).
\end{equation}
\end{Theorem}

Notice that $\forall x^\ast \in A^\ast$, there exists a sequence of partitioning sets $\{A_{i_k, k}: k=1,...\}$ s.t. $x^\ast \in A_{i_k, k}$ and $A_{i_k, k} \in \pi_k$, since the sets in each partition $\pi_k$ completely cover $A_k$ and hence $A^\ast$.
Furthermore, by the properties of the exponential family of densities (see (\ref{eq.3.3})), we have that
\begin{equation*}
\mathbb{E}_{h_{{i_k},k}}[\Gamma(x)]=\mathbb{E}_{\theta_{i_k, k}}[\Gamma(x)],
\end{equation*}
where recall that $h_{{i_k},k}\propto \mathds{1}\{x\in A_{{i_k},k}\}$.

Therefore, to show (\ref{eq.3.7}), it remains to show
\begin{equation*}
\lim_{k\rightarrow\infty}\mathbb{E}_{h_{i_k, k}}
\left[\Gamma(x)\right]=\Gamma(x^\ast),
~~~\text{or equivalently,} ~~~
\lim_{k\rightarrow\infty}\mathbb{E}_{h_{i_k, k}}
\left[\Gamma(x)-\Gamma(x^\ast)\right]=0.
\end{equation*}
That is, to show
\begin{equation}\label{eq.3.8}
\lim_{k\rightarrow\infty}\int_{A_{i_k, k}}
\left(\Gamma(x)-\Gamma(x^\ast)\right)h_{i_k, k}(x)dx=0.
\end{equation}
Given that $\Gamma(\cdot)$ is continuous on $\mathcal{X}$, we have that $\forall \epsilon>0$, $\exists \delta>0$ s.t.
\begin{equation*}
\abs{\Gamma(x)-\Gamma(x^\ast)}\le\epsilon, ~~\forall x\in B_{\delta}(x^\ast),
\end{equation*}
where $B_{\delta}(x^\ast) := \{x\in\mathcal{X}:
\norm{x-x^\ast}\le \delta\}$ represents the neighborhood ball centered at $x^\ast$ with radius $\delta$.
Further note that $\lim_{k\rightarrow \infty} \lvert \pi_k \rvert =0$, we have
\begin{equation*}
\lim_{k\rightarrow \infty} diam (A_{i_k, k}) =0.
\end{equation*}
Therefore, there exists a large integer $K_\epsilon$ depending on $\epsilon$ such that for all $k\ge K_\epsilon$, $A_{i_k, k}\subseteq B_{\delta}(x^\ast)$, we note that $x^\ast \in A_{i_k, k}$.
It follows that for all $k\ge K_\epsilon$,
\begin{equation*}
\int_{A_{i_k, k}}
\left(\Gamma(x)-\Gamma(x^\ast)\right)h_{i_k, k}(x)dx \le
\epsilon \int_{A_{i_k, k}}
h_{i_k, k}(x)dx =\epsilon.
\end{equation*}
Therefore, (\ref{eq.3.8}) and hence Theorem \ref{thm.3.1} holds.

Theorem \ref{thm.3.1} implies that for any solution with domination measure zero, there exists a sequence of exponential sampling distributions that converges to a degenerate distribution on that solution.

In practice, Algorithm ``$\mathbf{SASMO_0}$'' is not implementable for the following reasons: 1) the set $A_k$, which is regarded as the promising region of the solution space, could not be constructed explicitly since the domination measure $D(x)$ is unknown;
2) solving for the sampling parameters exactly through minimizing KL divergence is unlikely since the reference sampling distribution $h_{i,k}(\cdot)$ do not have an explicit characterization.

\subsection{Algorithm---An Implementable Version}
\label{Sec3.subsec2:Implementable}

To have an implementable version of Algorithm ``$\mathbf{SASMO_0}$'', we integrate a sampling step at each iteration, in which multiple i.i.d. candidate solutions are drawn according to a certain sampling distribution and the corresponding values of domination measure are evaluated or estimated.

There many reasons why the sampling step is integral in constructing an implementable version of Algorithm ``$\mathbf{SASMO_0}$''.
First, the domination measure $D(x)$ can be approximated for all the candidate solutions, which will then be used to determine the reference value $\gamma_k$ and characterize the set of promising solutions $A_k$.
In particular, suppose $\mathbf{g}_k(\cdot)$ is the sampling distribution at the $k^{th}$ iteration, and $N_k$ i.i.d. candidate solutions $\{x^1_k,...,x^{N_k}_k\}$ are drawn according to $\mathbf{g}_k(\cdot)$.
Then the domination measure for each candidate solution can be estimated by
\begin{equation*}
\widetilde{D}(x_k^i)=\frac{1}{N_k\cdot\nu(\mathcal{X})}
\sum_{j=1}^{N_k}\frac{1}{\mb{g}_k(x_k^j)}
\mathds{1}\left\{x_k^j\prec_d x_k^i\right\}
\end{equation*}
via the principle of importance sampling.
Second, the partition $\pi_k$ will be determined based on the evaluations of candidate solutions;
that is, the partitioning sets $\{A_{i_k, k}\}$ will be characterized by clusters of the elite candidate solutions from a certain clustering algorithm. Thus, the reference sampling distributions $\{h_{i_k, k}(\cdot)\}$ will be characterized by empirical distributions consisting of the elite candidate solutions in those clusters. In this case the magnitude of $\pi_k$ will be determined by the parameters of the clustering algorithm. Third, the sampling parameters $\{\theta_{i_k, k}\}$ will be solved by minimizing the KL divergence between the constructed empirical sampling distributions and the parameterized sampling distributions. To this end, let us describe the following Algorithm \ref{alg.3.2}, which is referred to as Algorithm ``$\mathbf{SASMO_1}$'', for simulation optimization of multi-objective problem (\ref{eq.2.1}).

In the initialization step (step 1) of Algorithm ``$\mb{SASMO_1}$'', a common choice of the parameterized family of densities is the exponential family.
The initial sampling parameter $\theta_{1,0}$ should be chosen in such a way that the resulted sampling distribution is close to the uniform distribution on $\mathcal{X}$, so that the entire solution space will be evenly explored in the early iterations.
For example, if the parameterized family of densities is the family of multi-variate Gaussian distributions, then $\theta_{1,0}$ is characterized by the mean vector $\mu_{1,0}$ and covariance matrix $\Sigma_{1,0}$. To enforce global exploration of the entire solution space, $\Sigma_{1,0}$ needs to be relatively large.
The mixed coefficient $\alpha$, the percent quantile $\rho$, as well as the sample size sequence $\{N_k\}$ will affect the robustness and convergence of the algorithm. We will discuss their selections in the Numerical Experiments Section.

\begin{algorithm}[htb]
\caption{\textbf{S}tochastic \textbf{A}daptive \textbf{S}earch for \textbf{M}ulti-objective \textbf{O}ptimization---\textbf{I}mplementable}\label{alg.3.2}
1. \textbf{Initialization}:  Choose a parameterized family of densities $\{g(x;\theta): \theta\in \Theta\}$ with initial parameter $\theta_{1,0}$ with $I_0=1$. Specify a mixing coefficient $\alpha\in (0,1)$, a percent quantile $\rho$, a sample size sequence $\{N_k\}$. Set $k=0$.
\\
2. \textbf{Sampling}:  Draw $N_k$ i.i.d. candidate solutions $\{x_k^i : i=1, 2,..., N_k\}$ according to $\mb{g}(\cdot;\theta_k)$, where
\begin{equation}\label{eq.3.9}
\mb{g}_k(\cdot)\overset{\triangle}=(1-\alpha)g_k(\cdot)+\alpha U_{\mathcal{X}}(\cdot)
\end{equation}
is a mixed sampling distribution,
$g_k(x)=\frac{1}{I_k}\sum_{i=1}^{I_k} g(x;\theta_{i, k})$, and $U_\mathcal{X}(\cdot)$ is the uniform distribution on $\mathcal{X}$.
\\
3. \textbf{Estimation}: For $i=1,...,N_k$, estimate the domination measure $D(x_k^i)$ at $x_k^i$  by
\begin{equation}\label{eq.3.10}
\widetilde{D}(x_k^i)=\frac{1}{N_k\cdot\nu(\mathcal{X})}
\sum_{j=1}^{N_k}\frac{1}{\mb{g}_k(x_k^j)}
\mathds{1}\left\{x_k^j\prec_d x_k^i\right\}.
\end{equation}
Sort $\{\widetilde{D}(x_k^i)\}$ in ascending order, denoted by $\widetilde{D}(x_1^{(1)})
\le \widetilde{D}(x_1^{(2)})\le \cdots \le \widetilde{D}(x_1^{(N_k)})$. Set the reference value $\widetilde{\gamma_k}$ to be the sample $\rho$-percent quantile $\widetilde{D}(x_1^{(\lceil \rho N_k\rceil)})$, i.e.,
$\widetilde{\gamma_k}=\widetilde{D}(x_1^{(\lceil \rho N_k\rceil)})$, where $\lceil \rho N_k\rceil$ is the smallest integer that is greater than or equal to $\rho N_k$.
\\
4. \textbf{Updating}: Construct the set of elite candidate solutions by $\widetilde{A}_k :=\{x_k^i: \widetilde{D}(x^i_k) \le \widetilde{\gamma_k}\}$. Using a clustering algorithm (e.g., Algorithm \ref{alg.3.3}) to cluster $\widetilde{A}_k$ into clusters
$\widetilde{\pi_k}:=\{\widetilde{A}_{1,k},...,\widetilde{A}_{I_k,k}\}$.
Update the parameter $\theta_{i,k}$ based on the set of elite candidate solutions $\widetilde{A}_{i,k}$ by solving
\begin{equation}\label{eq.3.11}
\theta_{i,k}\overset{\triangle}=\argmax_{\theta\in \Theta}\frac{1}{\lvert \widetilde{A}_{i,k} \rvert}
\sum_{x\in \widetilde{A}_{i,k}} \frac{\ln g(x;\theta)}{\mb{g}_k(x)}.
\end{equation}
5. \textbf{Stopping}. Check if some stopping criterion is satisfied. If yes, stop and return the means of the current parameterized sampling distributions; else, set $k:= k+1$ and go back to step 2.
\end{algorithm}

In the sampling step (step 2), note that in (\ref{eq.3.9}) the sampling distribution $\mathbf{g}_k(\cdot)$ is derived by mixing
the uniform distribution on $\mathcal{X}$ with the equal-weighted combination of the parameterized sampling distributions obtained from the previous iteration. The uniform distribution component in $\mathbf{g}_k(\cdot)$ helps maintain a global exploration of the entire solution space. Furthermore, since the candidate solutions drawn from $\mathbf{g}_k(\cdot)$ will be used to estimate the domination measure, the uniform distribution component controls the variance in the estimations. In particular, the variances will be bounded and the bounds only depend on the choice of $N_k$. Lastly, the choice of mixing coefficient $\alpha$ affects the robustness of the algorithm since it determines how much global exploration is achieved and the amount of variance reduction. A typical choice of $\alpha$ is $\alpha=0.1$.

In the estimation step (step 3), for the sake of convergence, the sample size $N_k$ is either set to be a large constant or in a way such that $N_{k+1}=\tau N_k$ for certain $\tau>1$.
The domination measure of the sampled candidate solutions are estimated via the principle of importance sampling, where the importance sampling distribution is also the mixed distribution $\mb{g}_k(\cdot)$.
As mentioned earlier, due to the uniform distribution component in the sampling distribution, the variances of the resulted estimators are bounded. The quantile level $\rho$ controls the reference value $\widetilde{\gamma}_k$, i.e., determines the number of elite candidate solutions that are used to update the sampling distribution in the next iteration, and the trade-offs between the exploitation of the neighborhood of current best solutions and the exploration of the entire solution space.
For example, when a smaller $\rho$ is used, less elite candidate solutions are used in updating the sampling distribution, which results in less exploration in the solution space.
Furthermore, note that the simple evolution of $\widetilde{\gamma}_k$ does not guarantee that it convergences to zero. However, it performs well in numerical tests and its performance is competitive to the more sophisticated methods of constructing $\widetilde{\gamma}_k$ in other model-based optimization methods, e.g., the one in SMRAS from \cite{hu2008model}.

In the updating step (step 4), the promising region $A_k$ is characterized by the set of elite candidate solutions $\widetilde{A}_k$, and the partition $\pi_k$ is characterized by the set of clusters
$\widetilde{\pi}_k$ yielded from applying a clustering algorithm on $\widetilde{A}_k$. Note that $\widetilde{\pi}_k$ also needs to satisfy
\begin{equation*}
\lim_{k\rightarrow \infty} \lvert \widetilde{\pi}_k \rvert=0, ~~\text{where}~~
\lvert\widetilde{\pi}_k \rvert
=\max_{1\le i\le I_k} diam(\widetilde{A}_{i,k}).
\end{equation*}
If a threshold-based clustering algorithm is used as the clustering algorithm, then the sequence of the threshold distances should decrease to zero. We defer the discussion of the details of this clustering algorithm (Algorithm \ref{alg.3.3}). The updating rule (\ref{eq.3.11}) on $\theta$ can be regarded as a sampled version of the updating rule in Algorithm (\ref{alg.3.1}), in which the reference sampling distribution $h_{i,k}(\cdot)$ is replaced by the corresponding empirical distribution. In particular, if the parameterized family of densities is the family of multivariate Gaussian family with $\theta=(\mu, \Sigma)$, then the mean $\mu_{i,k}$ and the covariance matrix $\Sigma_{i,k}$ computed from (\ref{eq.3.11}) is the sample mean and covariance matrix given by the candidate solutions in $\widetilde{A}_{i,k}$. We point out the resulted sampling distribution parameter could be further projected onto a properly constructed subset of $\Theta$ for numerical stability.

In the stopping step (step 5), a common stopping criterion is when the threshold distance in the clustering algorithm falls below a pre-specified threshold bound $\bar{\Delta}$ or a maximum number of iterations $t_{max}$ is reached.
The resulted means of the current parameterized sampling distributions will form a finite and approximate uniform representation of the solutions with domination measure of zero, with the hope that these solutions are close to the Pareto optimal set.
Moreover, the pre-specified threshold bound or maximum number of iterations will largely determine the number of solutions in the resulted approximation of the Pareto optimal set.

In a broader sense, we note that Algorithm ``$\mb{SASMO_1}$'' can be used to generate a finite and an approximate uniform representation of the optimal solution set for a single-objective optimization problem with multiple (possible uncountable) global optima, with a deterministic or stochastic objective.
Here the advantage of using Algorithm ``$\mb{SASMO_1}$'' for finding multiple minima of the stochastic minimization problem (\ref{eq.2.4}) lies in that only a single-layer of simulation is needed instead of a nested one, since the sampling is used to search the solution space and estimate the domination measure.
This is in contrast to using a model-based algorithm to solve a general stochastic optimization problem, because an outer-layer simulation is employed to sample on the solution space and an inner-layer simulation is needed for estimation of the function values.

The convergence analysis of Algorithm ``$\mb{SASMO_1}$'' is more complicated due to the involvement of sampling.
Strong requirements on the sample size sequence, the reference value sequence, the choice of clustering algorithm, and even the updating rule might need to be imposed to guarantee the convergence of the algorithm.
We leave the convergence analysis of  Algorithm ``$\mb{SASMO_1}$'' for future research.
Nevertheless, we will show that Algorithm ``$\mb{SASMO_1}$'' performs competitively to or outperforms some of the existing algorithms. Note that similar to Algorithm ``$\mb{SASMO_1}$'', these methods do not have convergence results.

We conclude this section by introducing the clustering algorithm used in the updating step of Algorithm ``$\mb{SASMO_1}$'', in which a threshold distance $\Delta_0$ and a threshold distance shrinking factor $C>1$ are chosen.
Of course, other clustering algorithms can be used as long as it satisfies the aforementioned guidelines.

\begin{algorithm}[htp]
\caption{A Threshold-Based Clustering Algorithm}\label{alg.3.3}
\textbf{Input}: Threshold distance $\Delta_k$ and elite solution set $\widetilde{A}_k$.
\\
\textbf{Output}: The set of clusters $\widetilde{\pi}_k$ and threshold distance $\Delta_{k+1}$.\\
1. \textbf{Initialization}: Randomly select a solution from $\widetilde{A}_k$. This solution is defined as the centroid of cluster $\widetilde{A}_{1,k}$.
\\
2. \textbf{Iteration}: Randomly select a solution from $\widetilde{A}_k$ that has not been assigned to any cluster.
Compute the Euclidean distances from the solution to the centroids of existing clusters in a randomized order.
Assign the solution to first cluster found with distance less than $\Delta_k$ and update the centroid of the cluster as the average of the solutions in the cluster. If no such cluster is found, create a new cluster where the solution is the centroid of the new cluster. \\
3. \textbf{Termination}: Check if there is solution from $\widetilde{A}_k$ that has not been clustered. If yes, go to step 2; otherwise return the set of clusters $\widetilde{\pi}_k$ and the threshold distance at next iteration by
\begin{equation}\label{eq.3.12}
\Delta_{k+1}=\min\left[\frac{1}{C I_k}\sum_{i=1}^{I_k}
Tr(\widetilde{\Sigma}_{i,k}),\frac{\Delta_k}{C}\right],
\end{equation}
where recall that $I_k$ is the number of clusters, $\widetilde{\Sigma}_{i,k}$ is the sample variance of cluster $\widetilde{A}_{i,k}$, and $Tr(\cdot)$ is the trace of a matrix.
\end{algorithm}

In the iteration step of Algorithm \ref{alg.3.3}, we randomized the order in which the distance from the selected solution and the centroids of the existing clusters is compared to prevent one cluster from becoming significantly larger than the other. In the termination step, the threshold distance at the next iteration $\Delta_{k+1}$ is decreasing adaptively, noting that $1/ I_k \cdot\sum_{i=1}^{I_k}Tr(\widetilde{\Sigma}_{i,k})$ is a empirical measure of how solutions within each cluster are close to each other.
If this measure is still large, then the threshold distance is forced to decrease by at least a factor of $C$. Therefore, the shrinking factor $C$ determines how fast Algorithm ``$\mb{SASMO_1}$'' terminates and how many solutions are generated in the final finite representation approximation of the Pareto optimal set.

\section{Numerical Experiments}\label{Sec4:Numerical}

The objective of our numerical results is to show that Algorithm ``$\mb{SASMO_1}$'' is
1) not sensitive to the geometry of the Pareto optimal set or the Pareto front,
2) scalable in terms of decision variables and objective functions, and 3) competitive to existing methods in terms of how close the solutions are to the true Pareto optimal set and how evenly spread the solutions are in the solution space.
To accomplish these goals, we evaluate the performance of Algorithm ``$\mb{SASMO_1}$'' on test functions from the ZDT (\cite{deb2001multi}), DTLZ (\cite{deb2001multi}), and Van Veldhuizen's (\cite{coello2002evolutionary}) test suites and compare our results with the following existing methods:

\begin{itemize}
\item Elitist Non-Dominated Sorting Genetic Algorithm (NSGA II) (\cite{deb2001multi})

\item Strength Pareto Evolutionary Algorithm (SPEA-II) (\cite{kim2004spea2})

\item Pareto Envelope-based Selection Algorithm (PSEA-II) (\cite{corne2001pesa})

\item Multi-objective Particle Swarm Optimization (MOPSO) (\cite{coello2002evolutionary})
\end{itemize}

The above methods are all evolutionary algorithms and they differ in how the population of candidate solutions are selected and maintained.
In NSGA-II the population of solutions is divided based on the following rule: the first group of solutions is all non-dominated and the second group of solutions is only dominated by the solutions in the first group.
The grouping is continued in this manner until all solutions are classified.
Once the solutions are divided, a crowding distance is calculated, which measures how close a solution is to its neighbors.
Solutions are selected based on their group classification and crowding distance, and new solutions are generated from crossover and mutation operators.
SPEA-II is an extended version of the original SPEA algorithm (\cite{zitzler1998}), which includes a specialized ranking system to order the solutions based on their fitness value, which is an objective function that summaries how close a given solution is to the Pareto front.
SPEA-II keeps an archive of all solutions generated starting from the initialization of the algorithm and constructs a population of solutions that combines the archived solutions with the solutions generated at the current iteration.
All non-dominated solutions in the population  are assigned fitness values such that the search is directed towards the true Pareto front.
PSEA-II introduces a new selection technique where the objective space is divided into hyperboxes and solutions are randomly selected from those hyperboxes.
The fitness value of a non-dominated solution depends on the number of non-dominated solutions that occupy that same hyperbox. This method of selection is shown to result in a good spread of solutions in the objective space.
MOPSO is a particle swarm method that includes a constraint-handling mechanism and a mutation operator that substantially improves the exploration ability of the original algorithm.

The problems in the ZDT test suite have a scalable number of parameters.
Therefore, this test suite tests the ability of an algorithm to converge to the Pareto front and obtain diverse solutions in a high-dimensional solution space.
The challenge in dealing with a high dimensional solution space is that it is more difficult to get an evenly spread of solutions in the solution space.
In practice the structure of the Pareto front is unknown; therefore, it is important to see how our method performs with different Pareto front geometries.
Consequently, we consider test function from the Van Veldhuizen's test suite since they present a variety of Pareto front geometries.
The different structures of the Pareto front can be convex, concave, degenerate, mixed, continuous, discontinuous, or contain flat regions. In this context, flat regions are areas of the objective space where relative small perturbations of parameters in the decision space do not affect the objective values.
Each geometric structure presents is own difficulty.
For instance, problems with isolated points are difficult to solve because usually there is no information in the surrounding region that indicates a Pareto optimal solution is nearby.
Furthermore, a function that has a many to one mapping between the decision space and the objective space are difficult to solve due to the flat regions in the objective space.
The last set of problems that we consider is the DLTZ test suite, which has a scalable number of objectives while also having complicated Pareto front geometries.
The increase in dimensionality of the objective space causes problems with selecting the best solutions.
When the number of objectives are large it causes a majority of solutions to be non-dominated by each other, which may
throttle an algorithm's convergence to the true Pareto front.
The problems chosen from the aforementioned test suites are included in the appendix and more details on the problems properties and challenges can be found in \cite{huband2006review}.

\begin{figure}\label{Graphs}
\centering
\begin{subfigure}{.5\textwidth}
  \centering
  \includegraphics[width=1\linewidth]{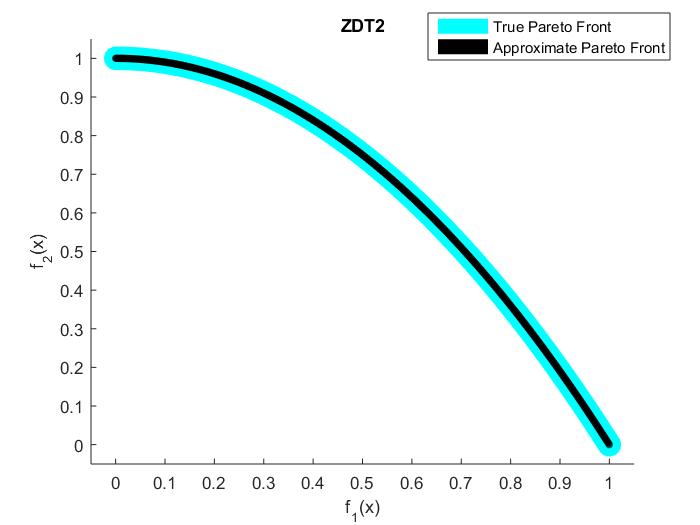}
  \caption{ZDT2}
  \label{fig.3.1}
\end{subfigure}%
\begin{subfigure}{.5\textwidth}
  \centering
  \includegraphics[width=1\linewidth]{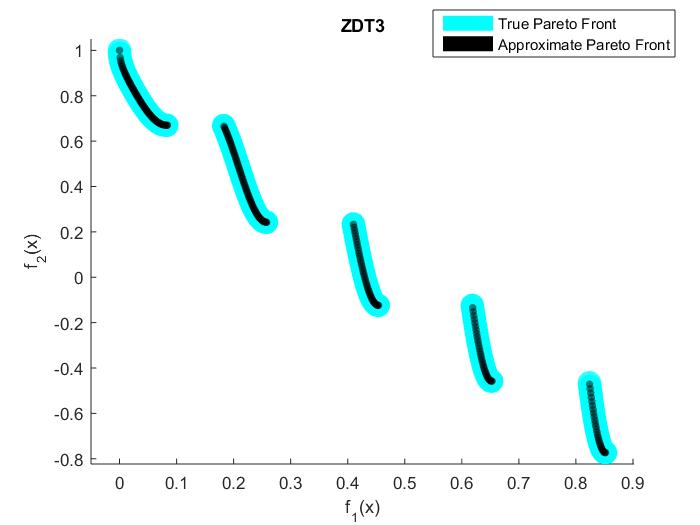}
  \caption{ZDT3}
  \label{fig.3.2}
\end{subfigure}
\begin{subfigure}{.5\textwidth}
  \centering
  \includegraphics[width=1\linewidth]{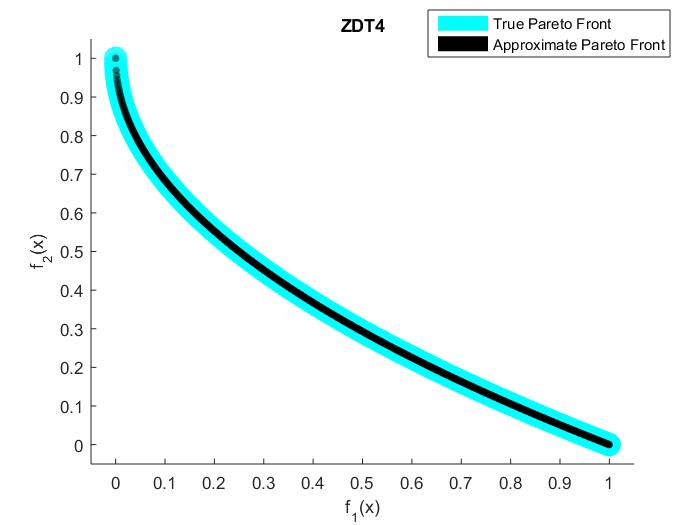}
  \caption{ZDT4}
  \label{fig.3.3}
\end{subfigure}%
\begin{subfigure}{.5\textwidth}
  \centering
  \includegraphics[width=1\linewidth]{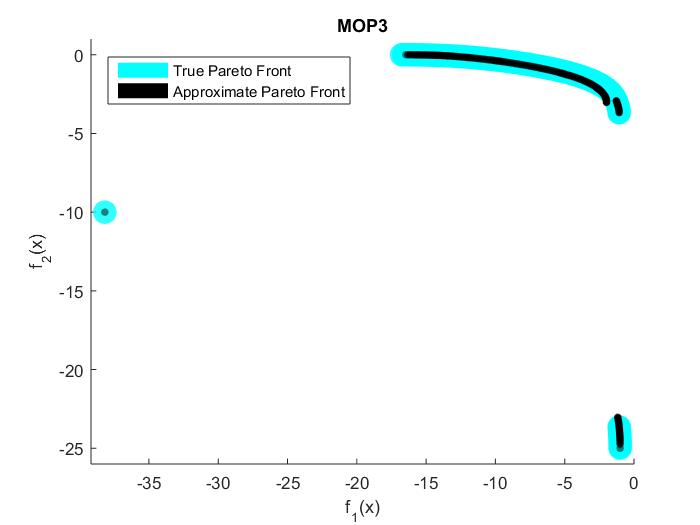}
  \caption{MOP3}
  \label{fig.3.4}
\end{subfigure}
\begin{subfigure}{.5\textwidth}
  \centering
  \includegraphics[width=1\linewidth]{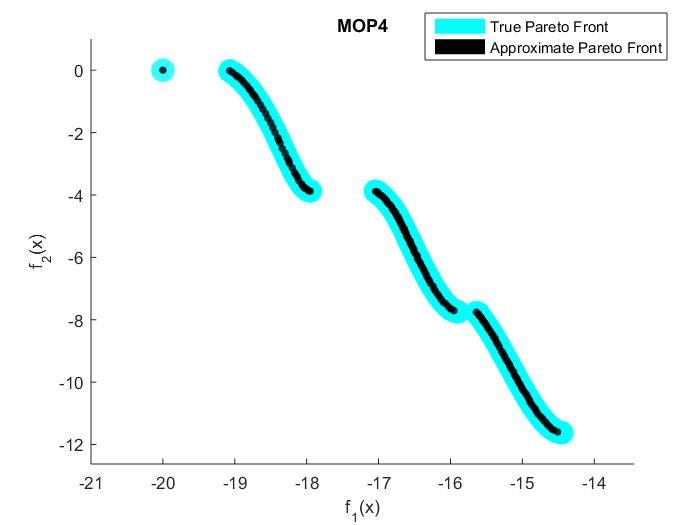}
  \caption{MOP4}
  \label{fig.3.5}
\end{subfigure}%
\begin{subfigure}{.5\textwidth}
  \centering
  \includegraphics[width=1\linewidth]{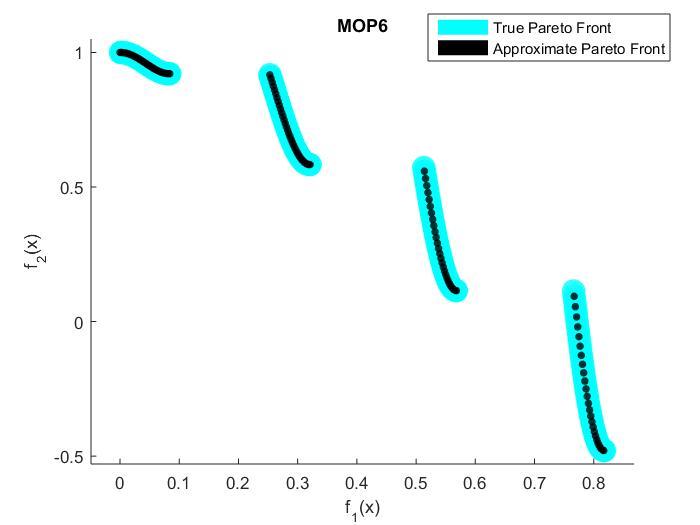}
  \caption{MOP6}
  \label{fig.3.6}
\end{subfigure}
\caption{Approximate Pareto Front V.S. True Pareto Front.}\label{fig.3.0}
\end{figure}

We solve all problems with Algorithm ``$\mb{SASMO_1}$'', which terminates when the threshold distance falls below the threshold bound $\bar{\Delta}$ or when the maximum number of iterations $t_{max}$ is reached.
We choose the following parameters: an initial sample size $N_0=1000$ and $N_k=k^{1.01} N_0$, a percent quantile $\rho=0.10$, a mixed coefficient $\alpha=0.1$, a threshold bound $\bar{\Delta} = 0.001$, a threshold distance shrinking factor $C =1.1$, an initial mean $\mu_0 = \mathbf{0}$, a maximum number of iterations $t_{max}=100$, and an initial covariance matrix $\Sigma_0=1000\mathbf{I}_d$.
The resulted approximations of the Pareto fronts for the tested problems are illustrated in Figure \ref{fig.3.0}.
Note that for each subfigure of Figure \ref{fig.3.0}, the cyan thick curve represents the true Pareto front and the black thin curve represents the approximated Pareto front produced by Algorithm ``$\mb{SASMO_1}$'' in the objective space $(f_1, f_2)$.
For the tested cases, we observe that Algorithm ``$\mb{SASMO_1}$'' is capable of obtaining isolated Pareto optimal points and capturing the entire Pareto front for problems that have multiple discontinuous Pareto curves.
These results also demonstrate that relaxing the concept of Pareto dominance to domination measure does not affect the solution quality of our algorithm.

In order to qualitatively quantify the performance of Algorithm ``$\mb{SASMO_1}$'' and compare it with the performances of the aforementioned existing algorithms, we use the convergence metric $\Lambda$ and diversity metric $\Upsilon$ defined in \cite{deb2001multi}.

The convergence metric $\Lambda$ in the objective space is defined by
\begin{equation*}
\Lambda \overset{\triangle}= \frac{1}{|R|} \sum_{y:~ \mb{f}(y)\in R}\left\{\min_{x \in Z}\norm{\mb{f}(x)-\mb{f}(y)}\right\},
\end{equation*}
where $R$ is a pre-specified reference set consisting of $|R|$ uniformly spread points from the true Pareto front and we choose $R=500$, $\mb{f}(\cdot)=(f_1(\cdot),...,f_n(\cdot))$ is the vector of objective functions, and $Z$ is the set of approximate Pareto front generated by the algorithm of interest.
In other words, the convergence metric $\Lambda$ can be regarded as the average distance from all points in the reference set to the approximate Pareto front, which measures the closeness of the approximate Pareto front to the true Pareto front.
Therefore, the smaller the value is for $\Lambda$, the closer the approximate Pareto front is to the true Pareto front.

Before introducing the diversity metric $\Upsilon$ in the solution space, let us first order the obtained $|Z|$ approximate Pareto optimal solutions $\{x^1,...,x^{|Z|}\}$ generated from an algorithm of interest by $\{x^{(1)},...,x^{(|Z|)}\}$ such that $x^{(1)}_1\le \cdots\le x^{(|Z|)}_1$.
That is, they are ordered by the values of their first components.
We also let $x^{(1)}$ and $x^{(|Z|)}$ be the left and right boundary points of the approximate Pareto optimal set, and let $x^l$ and $x^r$ be the left and right boundary of the true Pareto optimal set also in terms of the value of a solution's first component. The diversity metric $\Upsilon$ is defined by
\begin{equation*}
\Upsilon = \frac{d_l+d_r+\sum_{i=1}^{|Z|-1} |d_i - \bar{d}|}{d_l +d_r +(|Z|-1)\bar{d}},
\end{equation*}
where $d_l :=\|x^l-x^{(1)}\|$ ($d_i:=\|x^{(i+1)}-x^{(i)}\|$) is the distance between the left (right) boundary point of the true Pareto optimal set and the left (right) boundary point of the approximate Pareto optimal set, $d_i:=\|x^{(i+1)}-x^{(i)}\|$ is the distance between an approximate Pareto optimal solution and its closest neighbor, and $\bar{d}:=1/(|Z|-1)\sum_{i=1}^{|Z|-1}d_i$ is the average of these distances. Essentially, $\Upsilon$ measures how well the solutions are evenly spaced in the solution space. The smaller the value is for this metric, the closer the approximate Pareto optimal set is from being uniformly distributed.

For each problem instance, we perform 30 independent replications of each method implemented in MATLAB. The codes for the existing methods can be found at \url{http://yarpiz.com/category/multiobjective-optimization}. We report the best values for the convergence metric $\Lambda$ and the diversity metric $\Upsilon$ obtained out of the 30 trials for each method.
For the existing methods, we choose the following parameters:
\begin{itemize}
\item Number of generations: 100
\item Population size: 1000
\item Archive size: 1000.
\end{itemize}

The results are summarized in Table \ref{Table4} and Table \ref{Table5}.

\begin{table}[htb]
\center
\scalebox{1}{
\begin{tabular}{|c|ccccc|}
\toprule
\hline
    Problem & $\mb{SASMO_1}$ & NSGA II & SPEA II & MOPSO & PSEA II \\
    \hline
    & $\Lambda$ & $\Lambda$ & $\Lambda$ & $\Lambda$ & $\Lambda$ \\
      \midrule
      ZDT2  & \textbf{.0023}  & .1064 & .0765 & .2134 & .1863  \\
      ZDT3  & \textbf{.1145} & .5657 & .3312 & .1821 & .2765  \\
      ZDT4  & \textbf{4.321}  & 7.5432  & 9.8756 & 10.321 & 9.216  \\
      MOP3  & \textbf{.0221} & .0265 &.0456  & .8973 & .0821 \\
      MOP4  & .0987 & .0345  & .0673 & .1268 & .0531  \\
      MOP5  & \textbf{.0299} & .0321 & .0312 & .0589 & .0439  \\
      MOP6  & .0532  & .0582 & .0439 & 1.2956 & .9882 \\
      DTLZ1  & 2.9831  & 2.4531 & 2.8742 & 3.9814 & 3.7124  \\
      DTLZ2  & \textbf{2.4389}  & 2.9875 & 2.6192 & 4.8621 & 4.0921  \\
        \bottomrule
        \hline
\end{tabular}
}
\caption{Comparisons of Convergence Metric $\Lambda$.}
 \label{Table4}
\end{table}

\begin{table}[htb]
\center
\scalebox{1}{
\begin{tabular}{|c|ccccc|}
\toprule
\hline
    Problem &{$\mb{SASMO_1}$} & SPEA II & NSGA II & MOPSO & PESA II\\
    \hline
    & $\Upsilon$ & $\Upsilon$ & $\Upsilon$ & $\Upsilon$ & $\Upsilon$ \\
      \midrule
      ZDT2  & \textbf{.5641}   &.9987 &.6534 &.9874 & .8761  \\
      ZDT3  &1.0912  & .7654 &.7543 &.7469 &.9113  \\
      ZDT4  &.9792 &.9986 &.8921 & .9653 & 1.0320  \\
      MOP3  &\textbf{.4289} & .4467  & .8543 & .8856 & 1.2349  \\
      MOP4  &.\textbf{3531}  & .3876& .5789 & .6721 & .9177  \\
      MOP5  &.\textbf{3119}  &.3998 & .5431 &.6321 &.7643  \\
      MOP6  &\textbf{.3793}  & .4189 & .6842 & .5432 & .6631  \\
      DTLZ1  &1.1341  & .8743 & 1.2380 & .8734 & .9921  \\
      DTLZ2  &1.5939 & .9213 &1.8965 &1.9334 &1.8663  \\
        \bottomrule
         \hline
\end{tabular}
}
\caption{Comparisons of Diversity Metric $\Upsilon$.}
 \label{Table5}
\end{table}

We can see that Algorithm ``$\mb{SASMO_1}$'' outperforms the existing methods on over half of the problems in respect to both the convergence and diversity metrics.
The favorable results w.r.t. the diversity metric is likely due to the fact that the center of each cluster represents an estimated Pareto optimal solution and it is at least the threshold distance away from the closest cluster.
As a result, the distance between each estimated Pareto optimal solution is close to the threshold bound for most of the problems. We point out that as the dimension of the solution space increases so does the diversity metric for all the algorithms.
Another takeaway from the numerical results is how well our algorithm performed on problems that had discontinuous Pareto fronts. This is likely due to the fact that an adaptive number of components is used in the mixed sampling distribution so that each promising region of the solution space is thoroughly explored.
In conclusion, we show that Algorithm ``$\mb{SASMO_1}$'' gives satisfactory results regardless of the geometry of the Pareto front, and is competitive to several existing algorithms.

\section{Conclusions and Future Research}\label{Sec5:Conclusions}

In this paper, we introduce a novel performance metric called domination measure to measure the quality of a solution in a multi-objective problem.
Instead of solving the original problem, we propose a model-based approach to find a finite and approximately uniformly spread representation of the solutions with domination measure of zero, which is close to a finite and approximately uniformly spread representation of the Pareto optimal set.
We present an ideal version algorithm that has nice convergence properties, and an implementable version of the algorithm that has competitive numerical performances compared to many existing approaches.
More sophisticated approaches and algorithms based on domination measure can be incorporated to further improve the results on theoretically convergence and numerical performance.
Another future direction of research is to extend Algorithm ``$\mb{SASMO_1}$'' to the setting of stochastic multi-objective optimization, where the sample size allocations across the candidate solutions could to be studied.

\section*{Acknowledgements}

This work was supported by National Science Foundation under Grants CMMI-1413790 and CAREER CMMI-1453934, and Air Force Office of Scientific Research under Grant YIP FA-9550-14-1-0059.

\newpage
\bibliographystyle{ormsv080}
\bibliography{Zhou-Bibtex}

\appendix
\section{Functions Tested in Numerical Experiments}
\begin{itemize}
\item ZDT2
\begin{eqnarray*}
\min~~ \mb{f}(x) &=& (f_1(x),f_2(x))\\
f_1(x) &=& x_1\\
f_2(x) &=& g(x)[1 -(x_1/g(x))^2],~~\text{where}\\
g(x)&=& 1 + 9 (\sum^{30}_{i=2}x_i)/(30-1)]\\
x_i &\in & [0,1],~~ i = 1,\ldots, 30.
\end{eqnarray*}
\item ZDT3
\begin{eqnarray*}
\min~~ \mb{f}(x) &=& (f_1(x),f_2(x))\\
f_1(x) &=& x_1\\
f_2(x) &=& g(x)[1 -(x_1/g(x))^2 - \frac{x_1}{g(x)}\sin (10 \pi x_1)], ~~\text{where}\\
g(x)&=& 1 + 9 (\sum^{30}_{i=2}x_i)/(30-1)\\
x_i  &\in& [0,1],~~ i = 1,\ldots, 30.
\end{eqnarray*}
\item ZDT4
\begin{eqnarray*}
\min~~ \mb{f}(x) &=& (f_1(x),f_2(x))\\
f_1(x) &=& x_1\\
f_2(x) &=& g(x)(1- \sqrt{x_1/g(x)}), ~~\text{where}\\
g(x) &=& 1+10(10-1)+\sum_{i=2}^{10}(x_i^2-10\cos(4\pi x_i))\\
x_i  &\in & [-5,5],~~ i = 1,\ldots 10.
\end{eqnarray*}
\item MOP3
\begin{eqnarray*}
\max~~ \mb{f}(x) &=& (f_1(x),f_2(x))\\
f_1(x) &=& -[1+(A_1-B_1)^2 +(A_2-B_2)^2]\\
f_2(x) &=& -[(x_1+3)^2 +(x_2+1)^2], ~~\text{where}\\
A_1 &=& 0.5\sin 1 - \cos 1 + 2\sin 2 - 1.5\cos 2\\
A_2 &=& 1.5\sin 1 - \cos 1 + 2\sin 2 - 0.5\cos 2\\
B_1 &=& 0.5\sin x_1 - 2\cos x_1 + \sin x_2 - 1.5\cos x_2\\
B_2 &=& 1.5\sin x_1 - \cos x_1 + 2\sin x_2 -0.5\cos x_2\\
x_i  &\in & [-\pi,\pi], ~~ i= 1 \ldots 2.
\end{eqnarray*}
\item MOP4
\begin{eqnarray*}
\min~~ \mb{f}(x) &=& (f_1(x),f_2(x))\\
f_1(x)&=& \sum^{2}_{i=1} (-10 \exp(-0.2\sqrt{x_i^2+x_{i+1}^2}))\\
f_2(x) &=& \sum^3_{i=1}(|x_i|^{0.8} +5\sin(x_i)^3)\\
x_i &\in& [-5,5],~~ i=1, \ldots, 3.
\end{eqnarray*}
\item MOP5
\begin{eqnarray*}
\min~~ \mb{f}(x) &=& (f_1(x),f_2(x),f_3(x))\\
f_1(x) &=& 0.5(x_1^2+x_2^2)+\sin(x_1^2+x_2^2)\\
f_2(x) &=& \frac{(3x_1-2x_2+4)^2}{8} + \frac{(x_1-x_2+1)^2}{27} +15\\
f_3(x) &=& \frac{1}{x_1^2+x_2^2+1} -1.1 \exp (-x_1^2 -x_2^2)\\
x_1, x_2  &\in& [-30,30].
\end{eqnarray*}
\item MOP6
\begin{eqnarray*}
\min~~ \mb{f}(x) &=& (f_1(x),f_2(x))\\
f_1(x) &=& x_1\\
f_2(x) &=& (1+10x_2)[1-(\frac{x_1}{1+10x_2})^{2} - \frac{x_1}{1+10x_2}sin(8\pi x_1)]\\
x_1, x_2  &\in& [0,1].
\end{eqnarray*}
\item DTLZ1
\begin{eqnarray*}
\min~~ \mb{f}(x) &=& (f_1(x),f_2(x),f_3(x))\\
f_1(x) &=& 1/2 x_1x_2(1+g(x))\\
f_2(x) &=& 1/2 x_1(1-x_2)(1+g(x))\\
f_3(x) &=& 1/2(1-x_1)(1+g(x)), ~~\text{where}\\
g(x)&=& 100(\norm{x} -2 + \sum^7_{i=3}((x_i-0.5)^2-\cos(20\pi(x_i-0.5)))\\
x_i  &\in& [0,1]~~, i =1,\ldots ,7.
\end{eqnarray*}
\item DTLZ2
\begin{eqnarray*}
\min~~ \mb{f}(x) &=& (f_1(x),f_2(x),f_3(x))\\
f_1(x) &=& \cos(x_1 \pi /2)\cos(x_2 \pi /2)(1+g(x))\\
f_2(x) &=& \cos(x_1 \pi /2)\sin(x_2 \pi /2)(1+g(x))\\
f_3(x) &=& \sin(x_1 \pi /2)(1+g(x)), ~~\text{where} \\
g(x)&= &\sum^{12}_{i=3} (x_i - 0.5)^2\\
x_i  &\in& [0,1],~~ i = 1,\ldots, 12.
\end{eqnarray*}
\end{itemize}

%

\end{document}